\documentclass[12pt]{amsart}
\def\BibTeX{{\rm B\kern-.05em{\sc i\kern-.025em b}\kern-.08em
    T\kern-.1667em\lower.7ex\hbox{E}\kern-.125emX}}
\usepackage{amsmath}
\usepackage{enumerate}
\usepackage{amsthm}
\newtheorem{thm}{Theorem}

\usepackage{amssymb}
\usepackage{amsfonts}
\usepackage{latexsym}
\usepackage{array}
\usepackage{verbatim}
\usepackage{mathptmx} 
\usepackage{graphicx} 
\usepackage{color}
\usepackage{fancyhdr} 
\usepackage{hyperref}
\hypersetup{
    colorlinks = true,
    urlcolor = blue,
    citecolor	 = blue
}
\title{\sc Values of twisted Artin $L$-functions}
\author{Kenneth Ward}
\address{NYU-ECNU Institute of Mathematical Sciences\\
New York University Shanghai\\
Shanghai, China\\
200062}
\email{kaw354@nyu.edu}

\date{}
\allowdisplaybreaks
\begin{document}
\dedicatory{Dedicated to the memory of Paul J. Sally, Jr.}
\begin{abstract}
This note gives a simple proof that certain values of Artin's $L$-function, for a representation $\rho$ with character $\chi_\rho$, are stable under twisting by an even Dirichlet character $\chi$, up to the $\dim(\rho)$th power of the Gauss sum $\tau(\chi)$ and an element generated over $\mathbb{Q}$ by the values of $\chi$ and $\chi_\rho$. This extends a result due to J. Coates and S. Lichtenbaum.
\end{abstract}
\subjclass{11F67, 11F80, 11L05, 11M06}
\keywords{Artin $L$-function, twisted, Dirichlet character, Gauss sum, special values}
\maketitle
\small
\raggedbottom
\section{Introduction} We let $K|F$ denote a finite Galois extension of algebraic number fields, and $\chi$ will denote the character associated with a representation of $\text{Gal}(K|F)$. Via C.L. Siegel and H. Klingen \cite{Kli,Sie}, if $\chi$ denotes a Dirichlet character and $K$ is a totally real number field, then the Dirichlet $L$-function $L(\chi,s)$ lies in $\mathbb{Q}(\chi)$ when $s$ takes negative integer values. For a representation $\rho$, the associated character will simply be written as $\chi_\rho$. For general representations $\rho$ when $K$ is a finite Galois extension of $\mathbb{Q}$, J. Coates and S. Lichtenbaum \cite{CoLi} decomposed the factors at infinity in the functional equation \cite[VII.12.6]{Neu} \begin{align}\label{eq1} \Lambda(K|F,\chi_\rho,s) = W(\chi_{\rho})\Lambda(K|F,\overline{\chi}_\rho,1-s) && (W(\chi_\rho)| = 1)\end{align} to show that, at a negative integer $s = m$ that is a critical point for the Artin $L$-function $L(K|\mathbb{Q},\chi_\rho,s)$, either \begin{enumerate}[(I)] \item $m$ is odd, and the fixed field $K_\rho$ of the kernel of $\rho$ is totally real; or \item $K_\rho$ is totally imaginary, conjugation is central in $\text{Gal} (K_\rho|\mathbb{Q})$, and $\chi_\rho(\sigma) = -\dim(\rho)$. \end{enumerate} They then employed Brauer's theorem on induced characters and functorial properties of Artin's $L$-function to prove that $L(K|\mathbb{Q},\chi_\rho,m) \in \mathbb{Q}(\chi_\rho)$ when $m$ is a negative integer and critical point of $L(K|\mathbb{Q}, \chi_\rho, s)$ (for all other negative integers $m$, this $L$-function trivially takes the value zero).

We consider a finite Galois extension $K|\mathbb{Q}$, and a representation $\chi_\rho$ of $\text{Gal}(K|\mathbb{Q})$. We let $\chi$ denote an even Dirichlet character, which is viewed as either a one-dimensional character of $\text{Gal}(K|\mathbb{Q})$ acting trivially on the conjugation automorphism, or any one-dimensional character of $\text{Gal}(K|\mathbb{Q})$ if $K$ is totally real. We also let $\mathbb{Q}(\chi_\rho,\chi)$ denote the field generated over $\mathbb{Q}$ by the values of $\chi_\rho$ and $\chi$. The character $\chi \otimes \chi_\rho = \chi \cdot \chi_\rho$ is associated with the representation $\chi \otimes \rho = \chi \cdot \rho$. With $\mathfrak{f} = \mathfrak{f} (\chi)$ the Artin conductor and $K_\chi$ the fixed field of the kernel of $\chi$, we recognize the character $\chi$ as a character of the ray class group $J_\mathbb{Q}^{\mathfrak{f}}/P_\mathbb{Q}^{\mathfrak{f}}$ via the Artin symbol $(K_\chi|\mathbb{Q},\cdot): J_\mathbb{Q}^{\mathfrak{f}}/P_\mathbb{Q}^{\mathfrak{f}} \rightarrow \text{Gal}(K_\chi|\mathbb{Q})$. The Gauss sum of $\chi$ is defined in the usual way: For an algebraic number field $F$ with ring of integers $\mathfrak{o}_F$, a Dirichlet character $\chi$ of $F$ with finite part $\chi_f$ as a Gr\"{o}{\ss}encharakter and conductor $\mathfrak{f}$, $\mathfrak{d}_{F|\mathbb{Q}}$ the different of $F|\mathbb{Q}$, $y \in \mathfrak{f}^{-1} \mathfrak{d}_{F|\mathbb{Q}}^{-1}$, and $Tr_{F|\mathbb{Q}}(\cdot)$ the trace from $F$ to $\mathbb{Q}$, the Gauss sum of $\chi$ at $y$ is given by $$\tau(\chi,y) = \sum_{\substack{ x \text{ mod } \mathfrak{f} \\ (x,\mathfrak{f})=1}} \chi_f(x) e^{2\pi i Tr_{F|\mathbb{Q}}(xy)},$$ where the sum is taken over a set of representatives of $(\mathfrak{o}_F/\mathfrak{f})^*$. We also define $\chi$ to be zero for any ideals not relatively prime to $\mathfrak{f}$ and set $\tau(\chi) = \tau(\chi, 1)$. Finally, for a field $K$, a subfield $F \subset K$, and two elements $a,b \in K$, we will say that $a \sim_F b$ if $a = kb$ for some $k \in F$. We may now state our main result. \begin{thm} If $s = m$ is a negative integer and critical point of $L(K|\mathbb{Q},\chi_\rho,s)$, then $$L(K|\mathbb{Q},\chi \otimes \chi_\rho, 1 -m ) \sim_{\mathbb{Q}(\chi_\rho,\chi)} \tau(\chi)^{\dim(\rho)} L(K|\mathbb{Q},\chi_\rho,1-m).$$ \end{thm} Viewing $\chi_\rho$ as a character of $\text{Gal}(K_\rho|\mathbb{Q})$, we may write, by Brauer's theorem on induced characters \cite[5.2]{Col}, \begin{equation}\label{eq2} \chi_\rho = \sum_{i=1}^l n_i \text{Ind}_{H_i}^G(\chi_i), \end{equation} where $G = \text{Gal}(K_\rho|\mathbb{Q})$, and for each $i=1,...,l$, $\chi_i$ is a character of degree one of a subgroup $H_i \le G$ and $n_i \in \mathbb{Z}$. By Galois theory, for each $i=1,...,l$, let $F_i$ be the subfield of $K_\rho$ so that $H_i = \text{Gal}(K_\rho|F_i)$. We also let $G_i = \text{Gal}(K|F_i)$ for each $i=1,...,l$. The notation in the decomposition \eqref{eq2} will be used henceforth. We first relate the Gauss sums of restrictions of $\chi$ to $L$-function values. We will show that if $s = m$ is a negative integer and critical point of $L(K|\mathbb{Q},\chi_\rho,s)$, then \begin{equation} \label{eq122} \prod_{i=1}^l \left(\tau\left( \chi|_{G_i}\right)\right)^{n_i} \sim_{\mathbb{Q}(\chi_\rho,\chi)} \frac{L(K|\mathbb{Q},\overline{\chi}_\rho, 1 - m)}{L(K|\mathbb{Q},\overline{\chi \otimes \chi}_\rho, 1 - m)}.\end{equation} We will also show that the product of restricted Gauss sums is also related to the value of the unrestricted Gauss sum: \begin{equation} \label{eq123} \tau(\chi)^{\dim(\rho)} \sim_{\mathbb{Q}(\chi)} \prod_{i=1}^l \left(\tau\left(\chi|_{G_i}\right)\right)^{n_i} . \end{equation} It is no surprise that ramification groups have a role in this, as we must measure the action of characters on primes. The importance of requiring that $\chi$ be an even Dirichlet character is clear: this ensures that the critical points of $L(K|\mathbb{Q},\chi \otimes \chi_\rho,s)$ and $L(K|\mathbb{Q},\chi_\rho,s)$ are the same, by conditions (I) and (II).

\section{The Main Result} Our arguments take advantage of the interaction of Brauer's theorem with the Artin $L$-function as in the original proof of its meromorphicity. We shall also employ ideal numbers, Galois theory, and local behavior of characters. Let us turn immediately to the proof of Theorem 1.

 \begin{proof}[Proof of Theorem 1.] Let us fix $i \in \{1,...,l\}$. We may choose ideal numbers $d_i, f_i$ of $F_i$ so that the different $\mathfrak{d}_{F_i|\mathbb{Q}}$ of $F_i|\mathbb{Q}$ satisfies $\mathfrak{d}_{F_i|\mathbb{Q}} = (d_i)$ and the conductor $\mathfrak{f}_i = \mathfrak{f}(\chi_i)$ of $\chi_i$ satisfies $\mathfrak{f}_i = (f_i)$ \cite[VII.7]{Neu}. We let $p_i = (p_{i,\tau})_{\tau \in \text{Hom}(F_i,\mathbb{\mathbb{C}})}$ be the exponent of the infinite part of the decomposition of $\chi_i$ as a Gr\"{o}{\ss}encharakter, $$\chi_i((x)) = \chi_{i,f}(x) \chi_{i,\infty}(x) = \chi_{i,f}(x) N_{F_i|\mathbb{Q}}\left( \left( \frac{x}{|x|}\right)^{p_i}\right),$$ where $N_{F_i|\mathbb{Q}}(\cdot)$ denotes the norm from $F_i$ to $\mathbb{Q}$. As $\chi_i$ is a character of the group of ideals of $F_i$ via the Artin map, it may be viewed as a Dirichlet character, and thus has $p_{i,\tau} = 0$ unless $\tau \in \text{Hom}(F_i,\mathbb{R})$. The root number $W(\chi_i)$ from the functional equation \eqref{eq1} satisfies \begin{align} \notag W(\chi_i) &= \left[\omega_4^{Tr_{F_i|\mathbb{Q}}(\overline{p}_i)} N_{F_i|\mathbb{Q}}\left(\left( \frac{f_i d_i}{|f_i d_i |} \right)^{\overline{p}_i}\right)\right]^{-1} \frac{\tau(\chi_i)}{\sqrt{\mathcal{N}_{F_i|\mathbb{Q}}(\mathfrak{f}_i)}} \\& \notag = \pm \left[\omega_4^{Tr_{F_i|\mathbb{Q}}(\overline{p}_i)}\right]^{-1} \frac{\tau(\chi_i)}{\sqrt{\mathcal{N}_{F_i|\mathbb{Q}}(\mathfrak{f}_i)}},\end{align} where $Tr_{F_i|\mathbb{Q}}(\cdot)$ denotes the trace from $F_i$ to $\mathbb{Q}$, $\mathcal{N}_{F_i|\mathbb{Q}}(\cdot)$ the ideal norm from $F_i$ to $\mathbb{Q}$, and $\omega_4$ the complex number $i$, to avoid confusion with indices. We let $d_{F_i}$ denote the discriminant of $F_i$, and we let $$N_{i,+} = \frac{[F_i:\mathbb{Q}]}{2} + \frac{1}{2} \sum_{\substack{\mathfrak{p}|\infty \\ \mathfrak{p} \text{ real}}} \chi_i(\phi_{\mathfrak{P}}), \;\;\;\;\;\; N_{i,-} = \frac{[F_i:\mathbb{Q}]}{2} - \frac{1}{2} \sum_{\substack{\mathfrak{p}|\infty \\ \mathfrak{p} \text{ real}}} \chi_i(\phi_{\mathfrak{P}}),$$ where for each prime $\mathfrak{p}$ appearing in these sums, $\mathfrak{P}$ may be any choice of prime of $K$ so that $\mathfrak{P}|\mathfrak{p}$ and $\phi_{\mathfrak{P}}$ is a generator of the decomposition group $G_{\mathfrak{P}}(K|F_i)$. This yields the following expression for the functional equation: \begin{align}\label{eq3} L(K|F_i, \chi_i, 1 - s) = \pm & \left[\omega_4^{Tr_{F_i|\mathbb{Q}}(\overline{p}_i)}\right]^{-1} \frac{\tau(\chi_i)}{\sqrt{\mathcal{N}_{F_i|\mathbb{Q}}(\mathfrak{f}_i)}}[|d_{F_i}|\mathcal{N}_{F_i|\mathbb{Q}}(\mathfrak{f}_i)]^{s - \frac{1}{2}} \left(\cos\left(\frac{\pi s}{2}\right) \right)^{N_{i,+} }  \\& \notag \times \left(\sin\left(\frac{\pi s}{2}\right) \right)^{N_{i,-} } (2 (2\pi)^{-s} \Gamma(s) )^{[F_i:\mathbb{Q}]} L(K|F_i, \overline{\chi}_i, s). \end{align} All of this holds equally well with $\chi|_{G_i} \otimes \chi_i$ in place of $\chi_i$. For an automorphism $\sigma \in \text{Gal}(\overline{\mathbb{Q}}|\mathbb{Q}(\chi_\rho,\chi))$, induction by characters \eqref{eq2} and functoriality of the Artin $L$-function give \begin{equation} \label{eq4} \frac{L(K|\mathbb{Q}, \chi \otimes \chi_\rho, m)}{L(K|\mathbb{Q}, \chi_\rho, m)} = \alpha \left[\prod_{i=1}^l \left(\frac{\tau(\chi|_{G_i} \otimes \sigma(\chi_i))}{\tau(\sigma(\chi_i))}\right)^{n_i} \right]  \frac{L(K|\mathbb{Q}, \overline{\chi \otimes \chi_\rho}, 1-m)}{L(K|\mathbb{Q}, \overline{\chi_\rho}, 1-m)}, \end{equation} where $\alpha \in \mathbb{Q}$. The automorphism $\sigma$ holds the following fixed: \begin{enumerate}\item [(1.a)] the (rational-valued) infinite part of either Gr\"{o}{\ss}encharakter $\chi|_{G_i} \otimes \chi_i$ or $\chi_i$; \item [(1.b)] the conductors of $\chi|_{G_i} \otimes \chi_i$ or $\chi_i$; and \item [(1.c)] the action of $\chi|_{G_i} \otimes \chi_i$ or $\chi_i$ on the generators of decomposition groups for infinite primes. \end{enumerate} By the functional equation \eqref{eq3} without the use of $\sigma$, it follows from (1.a)-(1.c) that \begin{equation} \label{eq999} \frac{L(K|\mathbb{Q}, \chi \otimes \chi_\rho, m)}{L(K|\mathbb{Q}, \chi_\rho, m)}  = \alpha \left[\prod_{i=1}^l \left(\frac{\tau(\chi|_{G_i} \otimes \chi_i)}{\tau(\chi_i)}\right)^{n_i} \right]  \frac{L(K|\mathbb{Q}, \overline{\chi \otimes \chi_\rho}, 1-m)}{L(K|\mathbb{Q}, \overline{\chi_\rho}, 1-m)},\end{equation} so that \eqref{eq4} and \eqref{eq999} give \begin{equation} \label{eq5}\prod_{i=1}^l \left(\frac{\tau(\chi|_{G_i} \otimes \sigma(\chi_i))}{\tau(\sigma(\chi_i))}\right)^{n_i}  = \prod_{i=1}^l \left(\frac{\tau(\chi|_{G_i} \otimes \chi_i)}{\tau(\chi_i)}\right)^{n_i}.\end{equation} Furthermore, as in Theorem 1.2 of \cite{CoLi}, we have \begin{equation} \label{eq6} \frac{L(K|\mathbb{Q}, \chi \otimes \chi_\rho, m)}{L(K|\mathbb{Q}, \chi_\rho, m)} \in \mathbb{Q}(\chi_\rho,\chi),\end{equation} by none other than the fact that the value of each $L$-function in the numerator and denominator of \eqref{eq6} lies in $\mathbb{Q}(\chi_\rho,\chi)$. We let $\zeta$ be a root of unity, chosen so that $\mathbb{Q}(\zeta)$ contains all of the values of the characters $\chi$ and $\chi_i$, as well as any roots of unity appearing in the Gauss sums $\tau(\chi_i)$, $\tau(\chi|_{G_i} \otimes \chi_i)$, and $\tau(\chi|_{G_i})$, for all $i=1,...,l$. Let $\sigma \in \text{Gal}(\mathbb{Q}(\zeta)|\mathbb{Q}(\chi_\rho,\chi))$, so that $\sigma(\zeta) = \zeta^d$ for some integer $d$. We apply $\sigma$ to the expression \begin{equation} \label{eq121} \prod_{i=1}^l \left(\frac{\tau(\chi|_{G_i} \otimes \chi_i)}{\tau\left(\chi|_{G_i}\right)\tau(\chi_i)}\right)^{n_i}. \end{equation} We note that, for any Dirichlet character $\chi$ of an algebraic number field, $\tau(\chi,a) = \overline{\chi}_f(a)\tau(\chi)$ if $(a,\mathfrak{f}) = 1$, but also that $\tau(\chi,a) = 0 = \overline{\chi}_f(a)\tau(\chi)$ if $(a,\mathfrak{f})\neq 1$, as we have defined $\chi$ to be zero for ideals not relatively prime to $\mathfrak{f}$. Thus, for each $i=1,...,l$, \begin{enumerate} \item [(2.a)]  \begin{align*} \sigma(\tau(\chi|_{G_i} \otimes \chi_i)) &= \tau(\chi|_{G_i} \otimes \sigma(\chi_i),d)\\&= \left(\overline{\left(\chi|_{G_i} \otimes \sigma(\chi_i)\right)}_f(d)\right) \tau(\chi|_{G_i} \otimes \sigma(\chi_i));\end{align*} \item [(2.b)] $\sigma(\tau(\chi|_{G_i})) = \tau(\sigma(\chi|_{G_i}),d) = \tau(\chi|_{G_i},d)= \left(\overline{\left( \chi|_{G_i}\right)}_f(d)\right) \tau(\chi|_{G_i})$; and \item [(2.c)] $\sigma(\tau(\chi_i)) = \tau(\sigma(\chi_i),d) = \left(\overline{\sigma(\chi_i)}_f(d)\right) \tau(\sigma(\chi_i))$. \end{enumerate} This requires that $d$ be relatively prime to the conductors of each of these characters. By \eqref{eq5} and (2.a)-(2.c), we find that $\sigma$ fixes \eqref{eq121}. By \eqref{eq4}, \eqref{eq6}, we obtain \eqref{eq122}.

As in the previous argument, we again select an appropriate choice of root of unity $\zeta$, which in this instance contains the values of the characters $\chi|_{G_i}$ and $\chi$, as well as all roots of unity appearing in $\tau(\chi|_{G_i})$ and $\tau(\chi)$, and we consider some $\sigma \in \text{Gal}(\mathbb{Q}(\zeta)|\mathbb{Q}(\chi))$, where $\sigma(\zeta) = \zeta^d$ for some integer $d$. As before, $d$ must be relatively prime to the conductors of $\chi$ and $\chi|_{G_i}$, for all $i=1,...,l$. Similarly to previous arguments, we obtain \begin{equation} \label{eq11} \sigma\left(\frac{\prod_{i=1}^l (\tau(\chi|_{G_i}))^{n_i}}{\tau(\chi)^{\dim(\rho)}} \right) = \left(\frac{\prod_{i=1}^l \left(\overline{\left({\chi|_{G_i}}\right)}_f(d)\right)^{n_i}}{\overline{\chi}_f(d)^{\dim(\rho)}}\right) \left(\frac{\prod_{i=1}^l (\tau(\chi|_{G_i}))^{n_i}}{\tau(\chi)^{\dim(\rho)}} \right).\end{equation} Note that $\sigma$ makes no appearance in this final expression, as $\sigma$ fixes the values of $\chi$ and $\chi|_{G_i}$, for all $i=1,...,l$. Let $K_\chi$ and $K_{\chi|_{G_i}}$ denote the fixed fields of the kernels of $\chi$ and $\chi|_{G_i}$, for each $i=1,...,l$, respectively. We have for each $i=1,...,l$ that $(\chi|_{G_i})_f(d) = \chi|_{G_i}((d))$, where $(d)$ is viewed as an ideal of ${K_{\chi|_{G_i}}}$, and likewise that $\chi_f(d) = \chi((d))$, where $(d)$ is viewed as an ideal of ${K_\chi}$. We write $d = \prod_k p_k^{a_k}$ as the prime factorization of $d$ in $\mathbb{Z}$. As $d$ is relatively prime to the conductor of $\chi$, each prime $p_k$ is unramified in $K_\chi$. Thus, for a prime $\mathfrak{P}|\mathfrak{p}_k$ of $K$ with decomposition group $G_\mathfrak{P}$ and inertia group $I_{\mathfrak{P}}$ of $\mathfrak{P}|\mathfrak{p}_k$, the character $\chi$ may be viewed as a character of $G_\mathfrak{P} / I_\mathfrak{P}$. 

Let us fix some $i \in \{1,...,l\}$. For each $k$, we write $p_k = \prod_j {\mathfrak{P}_{k,j}}^{e_{\mathfrak{P}_{k,j}}}$ in $K$ and $p_k = \prod_m {\mathfrak{P}'_{k,m}}^{e_{\mathfrak{P}'_{k,m}}}$ in $F_i$. For each $j$, we write $\phi_{\mathfrak{P}_{k,j}}$ for a representative element of $G_{\mathfrak{P}_{k,j}}$ which maps to the Frobenius in the residue field of $\mathfrak{P}_{k,j}|\mathfrak{p}_k$, and we write $\phi'_{\mathfrak{P}_{k,j}}$ for the analogous representative in $\text{Gal}(K|F_i)$ \cite[VII.10]{Neu}. Let us also fix a choice of $k$ and $j$. Let $\mathfrak{P}'_{k,m}$ denote the prime of $F_i$ with $\mathfrak{P}_{k,j}|\mathfrak{P}'_{k,m}|\mathfrak{p}_k$, and let $f_{\mathfrak{P}'_{k,m}} = [\mathfrak{o}_{F_i}/\mathfrak{P}': \mathbb{Z}/p_k \mathbb{Z}]$. By definition of $\phi_{\mathfrak{P}_{k,j}}$ and $\phi'_{\mathfrak{P}_{k,j}}$, the image $$\overline{\phi_{\mathfrak{P}_{k,j}}^{f_{\mathfrak{P}'_{k,m}}}}(\alpha \text{ mod }\mathfrak{P}_{k,j})= \phi_{\mathfrak{P}_{k,j}}^{f_{\mathfrak{P}'_{k,m}}}(\alpha)\text{ mod }\mathfrak{P}_{k,j}$$ in the residue field of $\mathfrak{P}_{k,j}|\mathfrak{p}_k$ is equal to $ {\phi'_{\mathfrak{P}_{k,j}}}(\alpha)\text{ mod }\mathfrak{P}_{k,j}.$ Thus $\phi_{\mathfrak{P}_{k,j}}^{f_{\mathfrak{P}'_{k,m}}}$ and $\phi'_{\mathfrak{P}_{k,j}}$ belong to the same coset of $I_{\mathfrak{P}_{k,j}}$ in $G_{\mathfrak{P}_{k,j}}$. As $\chi$ is trivial on $I_{\mathfrak{P}_{k,j}}$, we have $\chi(\phi_{\mathfrak{P}_{k,j}}^{f_{\mathfrak{P}'_{k,m}}}) = \chi(\phi'_{\mathfrak{P}_{k,j}})$. By the identity $\sum_{m} e_{\mathfrak{P}'_{k,m}} f_{\mathfrak{P}'_{k,m}} = [F_i:\mathbb{Q}]$, we obtain $\chi((d))^{[F_i:\mathbb{Q}]} = \chi|_{G_i}((d))$. Returning to Brauer's theorem on induced characters \eqref{eq2}, we also have $\dim(\rho) = \sum_{i=1}^l n_i  [F_i: \mathbb{Q}]$. It follows that \begin{align} \label{eq17} {\prod_{i=1}^l \left({\left({\chi|_{G_i}}\right)}_f(d)\right)^{n_i}} = {{\chi}_f(d)^{\dim(\rho)}}. \end{align} By \eqref{eq11} and \eqref{eq17}, $\sigma$ fixes ${\prod_{i=1}^l (\tau(\chi|_{G_i}))^{n_i}}/{\tau(\chi)^{\dim(\rho)}}$, which must therefore lie in $\mathbb{Q}(\chi)$. This proves \eqref{eq123}. Finally, as $\tau(\chi) \tau(\overline{\chi}) \in \mathbb{Q}(\chi)$, we obtain from \eqref{eq122} and \eqref{eq123} that \begin{align*} L(K|\mathbb{Q},{\chi \otimes \chi}_\rho, 1 - m) \sim_{\mathbb{Q}(\chi_\rho,\chi)}  \tau(\chi)^{\dim(\rho)}  L(K|\mathbb{Q},{\chi}_\rho, 1 - m), \end{align*} proving Theorem 1. \end{proof}

\section{Acknowledgements} The author would like to thank New York University for its support while this paper was written, as well as Anantharam Raghuram for posing this question on Gauss sums and $L$-functions.

\end{document}